\documentclass[12pt,reqno]{amsart}
\usepackage{amsmath}
\usepackage{hyperref} 
\usepackage{amsfonts}
\usepackage{amssymb}
\usepackage{amsthm}
\usepackage[arrow, matrix, curve]{xy} 
\usepackage{hyperref}
\usepackage{svn-multi}
\usepackage{amssymb,amscd,url,tikz,stmaryrd,mathtools,fixltx2e}
\usepackage{color}
\usepackage{amsaddr}

\usepackage{paralist,csquotes,stackrel}
 
\usepackage[mathscr]{eucal}

\numberwithin{equation}{section}

\newtheorem{proposition}{\textbf{Proposition}}
\newtheorem{lemma}[proposition]{\textbf{Lemma}}

\newtheorem{theorem}[proposition]{\textbf{Theorem}}

\newtheorem{definition}[proposition]{\textbf{Definition}}

\newtheorem{remark}[proposition]{\textbf{Remark}}

\newcommand{\SL}{{\rm SL}}

\newcommand{\SU}{{\rm SU}}

\newcommand\C{{\mathbb C}}

\newcommand{\R}{{\mathbb R}}

\numberwithin{proposition}{section}

\numberwithin{equation}{section}

\begin{document}

\title[Projective structures and hyper-K\"ahler manifolds]{Real projective structures on Riemann surfaces  and new hyper-K\"ahler manifolds}

\author{Sebastian Heller}
\address{Institute of Differential Geometry,
Leibniz Universit\"at Hannover} 
 \email{seb.heller@gmail.com}


\date{\today}

\thanks{ We would like to thank J\"org Teschner for first pointing us to (integral) grafting of Fuchsian projective structures, and its interpretation as
constant curvature -1 metrics with singularities on compact Riemann surfaces.
We also thank the DFG for financial support through the research training group RTG 1670.}

\maketitle

\begin{abstract}
The twistor space of the moduli space of solutions of Hitchin's self-duality equations  can be identified with the Deligne-Hitchin moduli space of $\lambda$-connections. We use real projective structures on Riemann surfaces to  prove the existence of new components of real holomorphic sections of the Deligne-Hitchin moduli space. Applying the twistorial construction 
 we show the existence of new hyper-K\"ahler manifolds associated to any compact Riemann surface of genus $g\geq2$. These hyper-K\"ahler manifolds can be considered as moduli spaces of (certain) singular solutions of the self-duality equations.
\end{abstract}

\section{Introduction}
\label{intro}
Hyper-K\"ahler manifolds have been introduced by Calabi in the late 1970's. They are Riemannian manifolds whose holonomies are contained in $Sp(k),$ where $4k$ is the dimension of the manifold.  Compact examples of dimension 4 are classified (they are either a 4-torus or a $K3$-surface), and non-compact examples in dimension 4 are well-understood. Many examples of non-compact
hyper-K\"ahler manifolds arise as  moduli spaces of solutions to certain gauge theoretic equations. Most relevant for us
are the moduli spaces $\mathcal M_{SD}$ of solutions of Hitchin's self-duality equations \cite{HitchinSD} on a  compact Riemann surface.
The twistor space, the complex manifold parametrised over the 2-sphere of compatible complex structures, of $\mathcal M_{SD}$ has a complex analytic reincarnation as the Deligne-Hitchin moduli space of $\lambda$-connections on $\Sigma$, as was first pointed out by Deligne, see \cite{Simpson}. 
In general, the twistor space of a hyper-K\"ahler manifold makes it possible to reconstruct the hyper-K\"ahler manifold from complex analytic data, e.g., the underlying smooth manifold is given as a component of the space of real holomorphic sections of the twistor space \cite{HKLR}. Real sections in this component are called twistor lines.
Most important in the twistorial (re)construction of hyper-K\"ahler metrics is the special type of the (holomorphic) normal bundle of a twistor line: it is $\mathcal O(1)^{2k}\to\mathbb CP^1$, and locally, the space of real holomorphic section must be a real manifold of dimension $4k$ by Kodaira deformation theory. Moreover, evaluation at different $\lambda$ yield the complex structures and their K\"ahler forms, see \cite{HKLR}.

In the case of the self-duality moduli space, the twistorial construction not only parametrises the space of solutions.
Real sections in the Deligne-Hitchin moduli space also give rise to a complex-analytic construction of the 
solutions by loop group factorisation methods \cite[Theorem 3.6]{BHR}.
Simpson \cite{Simpson},\cite{Simpson2} has asked the natural question whether all real holomorphic sections are twistor lines, i.e., correspond to solutions of the self-duality equations. In \cite{HH2} and \cite{BHR} it was shown  
that there exists other real holomorphic sections besides the twistor lines. The class of counter-examples in \cite{HH2}
are given by solutions of the self-duality equations away from certain real curves on the given Riemann surface, while
the examples constructed in \cite{BHR} correspond to harmonic maps into the Lorentzian-symmetric deSitter 3-space
$\SL(2,\C)/\SU(1,1).$  The normal bundles of those real holomorphic sections have not been computed so far.

The aim of the paper is to show the existence of 
new hyper-K\"ahler manifolds associated to a Riemann surface of genus $g\geq2$.
We start with a real projective structure on the Riemann surface $\Sigma$, i.e., a complex projective structure with $P\SU(1,1)$-monodromy. We use the notation of \cite{Faltings}, and the reader should be aware that we do not consider $\R P^2$-structures.
All real projective structure are obtained by grafting \cite{Maskit,Hejhal,Sullivan-Thurstan}
with respect to a collection of disjoint simple closed curves up to isotopy \cite{Goldman}. The empty collection yields the self-duality moduli space.
We show that
a real projective structure gives rise to a section of the Deligne-Hitchin moduli space of $\Sigma$. This is analogous to uniformization of $\Sigma$,
which corresponds to the self-duality solution associated to the Higgs pair
\[(S\oplus S^*, \begin{pmatrix}0&0\\1&0\end{pmatrix}),\] where $S$ is a spin bundle on the Riemann surface.
The main difference to the case of uniformization is that the induced conformal metric of constant curvature -1 develops singularities (along curves where the developing map of the projective structure leaves the hyperbolic disc). In this regard, these real holomorphic sections are similar to those constructed in \cite{HH2}, i.e., they are self-duality solutions away from singularity curves on the Riemann surface. We next prove that the normal bundles of the sections obtained from real projective structures are
given by \[\mathcal O(1)^{6g-6}\to\C P^1,\]
where $g$ is the genus of $\Sigma.$ Applying the twistorial construction \cite{HKLR}
yields new hyper-K\"ahler metrics on spaces of real holomorphic sections. These hyper-K\"ahler manifolds can be considered as moduli spaces of (certain) singular solutions of the self-duality equations, see Theorem \ref{MT}.
One of the advantages of the  sections obtained from real projective structures
 (compared to those constructed in \cite{HH2})  is that their construction is quite simple and allows for explicit computations. Moreover, they provide a convenient conjectural picture of the space of components of real holomorphic sections of the Deligne-Hitchin moduli spaces, parametrised by isotopy classes of finite collections of disjoint simple non-trivial curves
 on a given Riemann surface, see Remark \ref{lastremark}.

\section{Background material}\label{sec:back}
In this section we shortly introduce the main objects of interest in this paper, namely Deligne-Hitchin moduli spaces and real projective structures obtained by grafting of Riemann surfaces.
\subsection{The Deligne-Hitchin moduli space}
We recall basic facts about $\lambda$-connections and Deligne-Hitchin moduli spaces. For more details we refer to
\cite{Simpson,Simpson2} where we first learnt about $\lambda$-connections, and to \cite{BHR,HH2}.
\begin{definition}
Let $\Sigma$ be a compact Riemann surface, and $V\to\Sigma$ be a complex vector bundle. A $\lambda$-connection
is a triple
\[(\lambda,\bar\partial,D)\]
where $\lambda\in\C,$ $\bar\partial$ is a holomorphic structure on $V$ and 
\[D\colon\Gamma(V)\to\Gamma(KV)\]
is a first order linear differential operator satisfying
\[D(fs)=\lambda\partial f \otimes s+ f Ds\]
for all functions $f$ and sections $s,$ and
\[\bar\partial D+D\bar\partial =0.\]
\end{definition}
A $\lambda$-connection $(\lambda,\bar\partial,D)$ for $\lambda\neq0$ gives rise to a flat connection
\[\bar\partial+\tfrac{1}{\lambda}D.\]
A $\lambda$-connection $(\lambda,\bar\partial,D)$ for $\lambda=0$ gives rise to a Higgs pair
\[(\bar\partial,\Phi=D).\]
A $\SL(2,\C)$ $\lambda$-connections is defined on a rank 2 vector bundle $V$ such that  the induced $\lambda$-connection on
$\Lambda^2 V$ is trivial. 
A $\SL(2,\C)$ $\lambda$-connection is called stable if there is no invariant line subbundle of non-negative degree and semi-stable if there is no
 invariant line subbundle of positive degree. A $\lambda$-connection with $\lambda\neq0$ is automatically semi-stable.
\begin{remark}
In this paper we only consider stable $\SL(2,\C)$ $\lambda$-connections, and we do not state that assumption 
explicitly in the following.  For all $\lambda$-connections which will be constructed in the paper, stability holds for obvious reasons or can be proven easily. We omit those proofs.
\end{remark}

We fix a topological trivialisation $V=\underline{\C^2},$ and consider the (complex) gauge group
\[\mathcal G^\C=\{g\colon \Sigma\to\SL(2,\C)\}.\]
The gauge group $\mathcal G^\C$ naturally acts on the space of $\lambda$-connections. Restricting to
stable $\lambda$-connections, the quotient is a complex manifold of dimension $6g-5$
\[\mathcal M_{Hod}=\mathcal M_{Hod}(\Sigma):=\{x=(\lambda,\bar\partial,D)\mid x \text{ is a stable } \lambda \text{ connection}\}/\mathcal G^\C,\]
see \cite{Simpson}. This space is called  the Hodge moduli space.
Elements in the Hodge moduli space are usually denoted by $[\lambda,\bar\partial,D]$ or by $[\lambda,\bar\partial,D]_\Sigma$
to emphasise their dependence on the Riemann surface $\Sigma.$
The Hodge moduli space admits a natural fibration to $\C.$

For a Riemann surface $\Sigma$ we denote by $\bar\Sigma$ the complex conjugate Riemann surface.
The Deligne gluing is a complex analytic diffeomorphism
\[\Psi\colon\mathcal M_{Hod}(\Sigma)_{\mid\C^*}\to \mathcal M_{Hod}(\bar\Sigma)_{\mid\C^*}\]
defined by
\[[\lambda,\bar\partial,D]_\Sigma\mapsto[\tfrac{1}{\lambda},\tfrac{1}{\lambda}D,\tfrac{1}{\lambda}\bar\partial]_{\bar\Sigma}.\]
The Deligne-Hitchin moduli space is defined as
\[\mathcal M_{DH}:=\mathcal M_{Hod}(\Sigma)\cup_\Psi\mathcal M_{Hod}(\bar\Sigma).\]
The Deligne-Hitchin moduli space admits a natural fibration to $\C P^1.$ The fiber over $\lambda=0$ is the (stable) Higgs bundle moduli space of $\Sigma$, the fiber over $\lambda=1$ is the moduli space of flat irreducible connections and the fiber over $\lambda=\infty$ is the moduli space of (stable) Higgs bundles on $\bar\Sigma.$
\subsubsection{Automorphisms of $\mathcal M_{DH}$} For every non-zero complex number $t\in\C^*$ there exists an automorphism
$t\colon \mathcal M_{DH}\to \mathcal M_{DH}$
determined by
\[t[\lambda,\bar\partial,D]_\Sigma:=[t\lambda,\bar\partial, t D]_\Sigma.\]
For $t=-1$ we denote this automorphism by $N.$

Denote the standard real structure on $V=\underline{\C^2}\to\Sigma$ by $\rho$. For a given linear differential operator
$A\in\mathcal D^k(\Sigma,V)$ we define
its complex conjugate 
\[\bar A:=\rho^{-1}\circ A\circ\rho=\rho\circ A\circ\rho\in\mathcal D^k(\Sigma,V).\]
With this terminology we obtain a real structure $\rho$ on $\mathcal M_{DH}$ defined by
\[\rho([\lambda,\bar\partial,D]_\Sigma)=[\bar\lambda,\overline{\bar\partial},\bar D]_{\bar\Sigma}.\]
The real involution covers
\[\lambda\mapsto 1/\bar\lambda.\]
Note that $N$ and $\rho$ commute, and hence give another real structure \begin{equation}\label{deftaus}\tau=N\circ\rho\colon\mathcal M_{DH}\to\mathcal M_{DH}\end{equation}
covering the antipodal involution $\lambda\mapsto -1/\bar\lambda$ of $\C P^1$. We remark that in higher rank $r>2$ there are different
possible definitions of the real involution $\rho$, see the discussion in \cite[$\S$ 1.6]{BHR}. 

\subsubsection{The twistor interpretation}
Hitchin's self-duality equations \cite{HitchinSD} on a Riemann surface are 
\[F^{\nabla+\Phi+\Phi^*}=0 \quad \text{ and } \quad \bar\partial^\nabla\Phi=0\]
for the $(0,1)$-part $\bar\partial^\nabla$ of a unitary connection $\nabla$ and $\Phi\in\Gamma(\Sigma,K End_0(V))$ and where $F$ denotes the curvature of a connection.
It  was shown by Hitchin that the moduli space $\mathcal M_{SD}$ of irreducible solutions to the self-duality equations modulo gauge transformations
is a hyper-K\"ahler manifold with respect to a natural $L^2$-metric.  The complex structures are induced by the Kobayashi-Hitchin correspondence: Every stable Higgs pair determines a unique solution (up to gauge transformations), giving the moduli space of solutions the complex structure $I$ of the Higgs bundle moduli space. On the other hand, every irreducible flat connection uniquely determines a solution of the self-duality equations \cite{Donaldson}, providing the complex structure $J$ on the moduli space of solutions. The complex structures $I$ and $J$ anti-commute, and $I$, $J$ and $K=IJ$ are K\"ahler with respect to a natural $L^2$-metric $G$.

The twistor space 
\[\mathcal P:= M\times \C P^1\to\C P^1\] of a hyper-K\"ahler manifold $(M,G,I,J,K)$ is naturally equipped with the almost complex structure 
\[I_{(p,\lambda)}:=(\frac{1-|\lambda|^2}{1+|\lambda|^2}I_p+\frac{\lambda+\bar\lambda}{1+|\lambda|^2}J_p+i\frac{\lambda-\bar\lambda}{1+|\lambda|^2}K_p,i).\]
This structure turns out to be integrable, see \cite{HKLR}. It admits a real involution
\[\tilde\tau\colon\mathcal P\to\mathcal P;\quad (p,\lambda)\mapsto (p,-\bar\lambda^{-1})\]
covering the antipodal involution. A holomorphic section $s$ of $\mathcal P\to\C P^1$ is called $\tilde\tau$-real if 
\[\tilde\tau (s(\lambda))=s(-\bar\lambda^{-1})\]
for all $\lambda.$
The manifold $M$ 
can be recovered as a component of the space of $\tilde\tau$-real holomorphic sections
of the holomorphic fibration
$\mathcal P\to\C P^1.$
In fact, for every $p\in M$
\[s_p\colon\C P^1\to\mathcal P;\quad \lambda\mapsto (p,\lambda)\]
provides a $\tilde\tau$-real holomorphic section, and locally, there cannot be any other $\tilde\tau$-real holomorphic section \cite{HKLR}.
Moreover, the Riemannian metric $G$ can be recovered from the twistor space. Consider the K\"ahler forms
$\omega_I$,  $\omega_J$ and $\omega_K$ with respect to $I,J$ and $K,$ respectively. Then,
\[\omega_\lambda:=(\omega_J+i\omega_K)+ 2\lambda \omega_I-\lambda^2(\omega_J-i\omega_K)\]
defines a holomorphic section
\[\omega_\lambda\in H^0(\mathcal P; \Lambda^{2,0} T^*F\otimes \mathcal O(2)),\]
where $\mathcal O(2)$ denotes the pullback of $\mathcal O(2)\to\C P^1,$ and $TF\to\mathcal P$ is the vertical bundle of
the fibration $\mathcal P\to\C P^1.$ The Riemannian metric $G$ is obtained  from $\omega_\lambda$ by evaluating derivatives of families of sections \cite{HKLR}. 
\begin{remark}
For $M=\mathcal M_{SD}$, the moduli space of solutions of the self-duality equations, $\omega_{\lambda=0}$ is the natural  symplectic form on the moduli space of Higgs bundles, while 
$\omega_{\lambda=1}$ is the Goldman symplectic form on the deRham moduli space of flat connections.
\end{remark}

A solution $(\nabla,\Phi)$ of the self-duality equations gives rise to its associated family of flat
connections \[\lambda\in\C^*\mapsto\nabla+\lambda^{-1}\Phi+\lambda \Phi^*,\] or likewise, to the holomorphic section 
\[\lambda\in\C\mapsto [\lambda,\bar\partial^\nabla+\lambda\Phi^*,\lambda\partial^\nabla+\Phi]_\Sigma\]
of the Deligne-Hitchin moduli space. In turns out that these sections are $\tau$-real \cite{Simpson}, see also \cite{BHR}.
\begin{theorem}[Deligne, Simpson, \cite{Simpson}]
The twistor space $\mathcal P\to\C P^1$ of the moduli space $\mathcal M_{SD}$ of solutions to the self-duality equations
on a compact Riemann surface $\Sigma$
is naturally biholomorphic to the Deligne-Hitchin moduli space via
\[(\nabla,\Phi,\lambda)\in \mathcal P=\mathcal M_{SD}\times\C P^1\mapsto [\lambda,\bar\partial^\nabla+\lambda \Phi^*,\lambda\partial^\nabla+\Phi]_\Sigma\in\mathcal M_{Hod}\subset\mathcal M_{DH}\]
such that  $\tilde\tau$ and $\tau$ coincide.
\end{theorem}
\subsubsection{Invariants of sections of the Deligne-Hitchin moduli spaces}
In \cite{BHR} several invariants for holomorphic sections $s\colon\C P^1\to\mathcal M_{DH}$ have been defined.
In particular, a section $s$ is called admissible if it admits a lift
\[\lambda\mapsto\nabla+\lambda^{-1}\Phi+\lambda\Psi\]
to a $\C^*$-family of flat connections
such that $(\bar\partial^\nabla,\Phi)$ and $(\partial^\nabla,\Psi)$ are stable Higgs pairs on $\Sigma$ and $\bar\Sigma$, respectively. 
Note that any section $s$ admits a lift of the form
\[\lambda\in\C^*\mapsto \lambda^{-1}\Phi+\nabla+\lambda\Psi_1+\lambda^2\Psi_2+...\]
with a stable Higgs pair $(\bar\partial^\nabla,\Phi)$. Also note that $\Psi_k,$ $k\in\mathbb N,$  are endomorphism-valued 1-forms and not necessarily of type $(1,0).$
The lift of a $\tau$-real holomorphic section yields a (holomorphic) family of $\SL(2,\C)$-valued gauge transformations
\[\lambda\in\C^*\mapsto g(\lambda)\in\mathcal G^\C\]
satisfying
\begin{equation}\label{eq:sign}\overline{\nabla^{-\bar\lambda^{-1}}}=\nabla^\lambda.g(\lambda)\end{equation} as a consequence of the reality condition.
By irreducibility of the connections \begin{equation}\label{eq:sign2}\overline{g(-\bar\lambda^{-1})}g(\lambda)=\pm \text{Id},\end{equation}
and the sign $\pm$ is an invariant of the section $s$. Note that the sign might change if we would allow $\mathrm{GL}(2,\C)$-valued gauge transformations.
A $\tau$-real holomorphic section is called positive or negative depending on the sign in \eqref{eq:sign2}.
Twistor lines in the Deligne-Hitchin moduli space are admissible negative sections.
In \cite{BHR} admissible positive $\tau$-real holomorphic sections of Deligne-Hitchin moduli spaces have been constructed, and 
in \cite{HH2} non-admissible  negative  $\tau$-real holomorphic sections of  Deligne-Hitchin moduli spaces have been constructed. The following theorem holds:
\begin{theorem} \cite{BHR}
An admissible $\tau$-real holomorphic section of a Deligne-Hitchin moduli space is either a twistor line or positive.
\end{theorem}

\subsection{Real projective structures}
Projective structures on Riemann surfaces, Schwarzian derivatives and opers are classical in the theory of Riemann surfaces. We shortly recall some basic facts, mainly to fix notations. For more details the reader might consult
\cite{Faltings,Goldman,BuPP}.

A projective structure on a Riemann surface is given by an atlas $(U_\alpha,z_\alpha)_{\alpha\in\mathcal U}$ whose transition functions are given by Moebius transformations, i.e., they satisfy
\[z_\beta\circ z_\alpha^{-1}(z)= \frac{a z+b}{c z+d}\]
for some (constant) \[\begin{pmatrix}a&b\\c&d\end{pmatrix}\in\text{GL}(2,\C)\] (depending on $\alpha,\beta\in\mathcal U)$.
While we stick to the notations of \cite{Faltings}, in \cite{Goldman} a projective structure is called a complex projective structure or a $\C P^1$-structure.

Obviously, $\C P^1$ is equipped with a projective structure, and a surface of genus $1$ is equipped with a complex structure via flat conformal coordinates.
A natural projective structure on any compact Riemann surface of genus $g\geq2$ is given by uniformization: The constant
curvature $-1$ metric provides a developing map to the hyperbolic disc $\mathbb H^2\subset\C P^1$ which is equivariant with respect to $P\SU(1,1)$-valued Moebius transformations.

 A projective structure on a compact Riemann surface can be described by a flat $\SL(2,\C)$-connection of a special form. Let $\nabla$ be a flat $\SL(2,\C)$-connection on $V=\underline\C^2\to\Sigma$ such that its induced holomorphic structure 
$\bar\partial^\nabla$ admits a holomorphic sub-line bundle $S$ of maximal possible degree $(g-1),$ $g$ being the genus of $\Sigma$.
Take a complementary $C^\infty$-bundle $S^*\subset V$, and decompose
\begin{equation}\label{eq:oper}\nabla=\begin{pmatrix} \nabla^{S}& \psi\\ \varphi &\nabla^{S^*}\end{pmatrix}\end{equation}
with respect to $V=S\oplus S^{*}.$ As $S$ is a holomorphic subbundle $\varphi$ is a $(1,0)$-form. Flatness implies that
\[\varphi\in H^0(\Sigma, K(S^{*})^2)\]
and $\varphi\neq0$ if $g\geq2$ because $S$ has positive degree. Hence,  $S$ is a spin bundle, i.e., $S^2=K$ as a holomorphic line bundle, and $\varphi$ is a nowhere vanishing section which is identified with a constant. The choice of the spin bundle $S$ corresponds to the choice of a lift of the 
$P\SL(2,\C)$-representation  to a $\SL(2,\C)$-representation.
\begin{definition}
A flat $\SL(2,\C)$-connection of the form \eqref{eq:oper} on a compact Riemann surface is called oper.
\end{definition}
The projective atlas corresponding to an oper $\nabla$ is obtained as follows: Consider two linear independent parallel sections
of $\nabla$
\[\Psi_1=\begin{pmatrix} x_1\\ y_1\end{pmatrix},\quad \Psi_2=\begin{pmatrix} x_2\\ y_2\end{pmatrix}\]
on an open set $U\subset \Sigma.$ Consider the projections of the sections to $V/S=S^*$. 
They are  holomorphic sections of $S^*$ as $V\to V/S$ is holomorphic. With respect to the decomposition $V=S\oplus S^*$ they are given as
$y_1$ and $y_2$. The quotient $z=y_1/y_2$ is a holomorphic map to $\C P^1$, and taking different parallel sections
\[\tilde\Psi_1=a\Psi_1+b\Psi_2 \quad  \tilde\Psi_2=c\Psi_1+d\Psi_2\]
(with $ad-bc=1)$ amounts into
\[z=y_1/y_2\mapsto \frac{a y_1+ b y_2}{c y_1+d y_2}=\frac{ az +b}{cz+d},\]
a Moebius transformation. Because $\varphi$  is nowhere vanishing $z$ is not branched, i.e., $z$ is a local (holomorphic) diffeomorphism. We obtain a projective atlas.

\subsubsection{Grafting} Grafting of projective structures was introduced by Maskit \cite{Maskit}, Hejhal \cite{Hejhal} and Sullivan-Thurstan \cite{Sullivan-Thurstan}. Our short description follows Goldman \cite{Goldman}.
Grafting yields projective structures whose monodromy is in $P\mathrm{SL}(2,\R)\cong \text{P}\SU(1,1)$. Following 
Faltings \cite{Faltings} and Takhtajan \cite{Tak} they are called {\em real projective structures} here. In \cite{Goldman}, a real projective structure is a $\R P(2)$-structure on the surface. The reader should be aware of the distinction.

For the construction of a real projective structure we consider the Fuchsian representation of a Riemann surface  $\tilde\Sigma$ and fix its developing map to $\mathbb H^2\subset\C P^1$. Consider
a simple closed curve in $\tilde\Sigma$ which is not null-homotopic. In its homotopy class there is a unique (oriented) geodesic $\gamma$ (with respect to the constant curvature -1 metric). Via the developing map $\gamma$ is mapped to a part of a circle.
The corresponding circle $C$ intersects the boundary at infinity of the hyperbolic disc at two points, one being the starting point $S$ and the other being the end point $E$ with respect to the orientation of the geodesic. The monodromy of the Fuchsian representation along $\gamma$ is given
by an element $A\in \SU(1,1)$, unique up to sign, for which $E$ is an attractive fixed point and $S$ is a repelling fixed point.
The sign depends on the lift of the Fuchsian representation to a $\SL(2,\C)$-representation. 
\begin{figure}
 \includegraphics[width=1\textwidth]{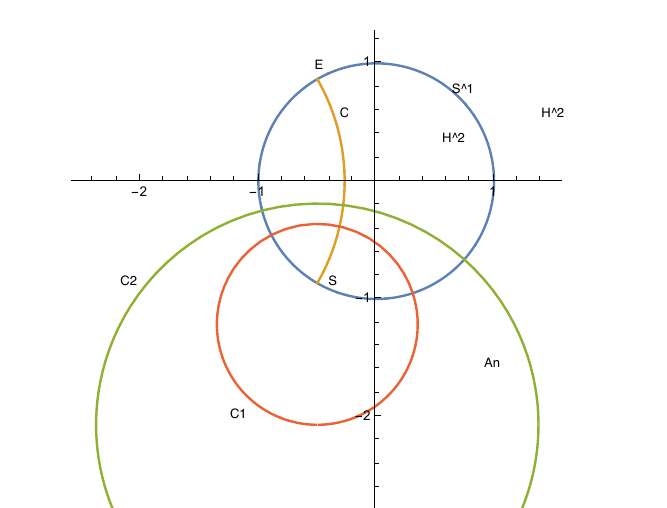}
 \caption{\footnotesize The annulus obtained from the Moebius transformation $A$.}\label{fig:an}
 \end{figure}

The transformation $A\in \SU(1,1)$ yields a torus $T$ with a flat  projective structure as follows:
Consider a circle $C_1$ which intersects both circles, $C$ and the boundary at infinity $S^1$, perpendicularly,
see Figure \ref{fig:an}.
It is mapped by $A$ to a circle $C_2=A(C_1).$ The circles $C_1$ and $C_2$ bound an annulus $An,$ and gluing the boundary circles
$C_1$ and $C_2$ via $A$ gives a torus $T$ with a complex projective structure. A fundamental piece of the developing map of this projective structure is the annulus $An.$

The torus $T$ can be glued to $\tilde\Sigma$ as follows: the geodesic arc of   $C$ contained in $\mathbb H^2$  intersects the annulus $An$  in a curve which projects to a closed curve $\tilde\gamma$ on $T$.  The torus $T$ admits
a metric of constant negative curvature $-1$ away from the intersection of its developing map with the boundary at infinity of $\mathbb H^2$, and $\tilde{\gamma}$ is a geodesic.
Denote by $|T|$ the annulus with two boundary components obtained by gluing two copies of $\tilde\gamma$ to the open annulus $T\setminus\tilde\gamma$.
There are open tubular neighbourhoods of $\gamma$ in $\tilde \Sigma$ and $\tilde\gamma$ in $T$ which are isometric. 
This isometry can be used to glue  the torus and $\tilde \Sigma$ to obtain a closed Riemann surface $\Sigma$ of the same genus by replacing the set $\gamma\subset \tilde\Sigma$ by $|T|$. Moreover,
$\Sigma$ has a projective structure induced from $\tilde\Sigma$ and $T$. This projective structure has the same monodromy as the original projective structure on $\tilde\Sigma.$ In particular, the monodromy is $\text{P}\SU(1,1)$-valued, and we obtain a real projective structure.
Note that (the developing map of) the projective structure on $\Sigma$ induces a curvature -1 metric away from a singularity set. The singularity set is given by two distinct  closed smooth
curves both homotopic to $\gamma$ (considered as a curve in $\Sigma$).

The above procedure is called grafting along $\gamma.$ We can  glue the torus multiple times, yielding yet another Riemann surface equipped with a projective structure. The monodromy is  still the Fuchsian monodromy of the initial Riemann surface $\tilde\Sigma$. Similarly, we can apply grafting to
the new Riemann surface $\Sigma$  along a  disjoint geodesic (with respect to the -1 curvature metric, disjoint from the singularity set).  Altogether, starting with a isotopy class $\mathcal C$ of simple not null-homotopic disjoint curves on $\tilde \Sigma$, and applying the grafting construction repetitively to the closed curves 
yields a new Riemann surface with a projective structure with $\SU(1,1)$-monodromy.

\begin{remark}\label{rem:sing}
Note that the topological information of $\mathcal C$ is still transparent after the grafting. In fact, the proof of the main theorem in \cite{Goldman} uses the decomposition of $\Sigma$ into connected components on which the induced Riemannian metric of constant curvature -1 is non-singular. The components are separated by closed curves on $\Sigma$ which
yield the information of $\mathcal C.$ 
\end{remark}
\begin{remark}
It was shown by Tanigawa \cite{Tanigawa} that grafting (with respect to an isotopy class $\mathcal C$ of simple not null-homotopic disjoint curves on the underlying topological surface) yields a homeomorphism of the Teichm\"uller space to itself. In particular, on any given Riemann surface $\Sigma$, there exists infinitely many projective structures with monodromy in $\text{P}\SU(1,1)$. The different real projective structures on a given Riemann surface $\Sigma$ are uniquely determined
by the choice of topological data $\mathcal C,$ i.e., by Goldman's result \cite{Goldman} every real projective structure arises from grafting in this way, and the isotopy class $\mathcal C$ can be recovered from the real projective structure from the set  of points where the (equivariant) developing map  of the projective structure intersects the boundary $\mathbb S^1=\partial\mathbb H^2\subset\mathbb CP^1$ at infinity.
\end{remark}
\begin{remark}
Real projective structures on a compact Riemann surface have been recently  identified \cite{Tesch} with eigenfunctions of the quantised Hamiltonians of the $\SL(2,\C)$ Hitchin system.
\end{remark}

\section{Real holomorphic sections via real projective structures}\label{sec:sections}
In this section we explain how real projective structures on $\Sigma$ give rise to  real holomorphic sections of the Deligne-Hitchin moduli space $\mathcal M_{DH}(\Sigma)$. We then
figure out some important properties of these sections.

Let $\Sigma$ be a compact Riemann surface of genus $g$ equipped with a real projective structure obtained by grafting along $\mathcal C.$ After the choice of a spin bundle $S$, the projective structure is given by a 
 flat $\SL(2,\C)$-connection $\nabla$ on $\Sigma$ with $\SU(1,1)$ monodromy. The underlying holomorphic structure
$\bar\partial^\nabla$ admits $S$ as a holomorphic subbundle.
The oper $\nabla$ differs from the uniformization connection $\nabla^{Fuchs}$ (for which the lift from $P\SL(2,\C)$ to $\SL(2,\C)$ is   determined by the same spin bundle $S$) by a holomorphic quadratic differential $q\in H^0(\Sigma,K^2).$
With respect to the induced $\SU(1,1)$-structure of the uniformization and the induced decomposition $V=S\oplus S^*$ into (complementary) orthogonal subbundles the uniformization oper is
\begin{equation}\label{nablafudec}
\nabla^{Fuchs}=\begin{pmatrix} \nabla^S& \varphi^*\\ \varphi &  \nabla^{S^*}\end{pmatrix}\end{equation}
where 
\[\varphi=1\in H^0(\Sigma,KK^{-1}),\]
$\varphi^*$ is the adjoint with respect to the metric with constant curvature $-1$.
Then, $\nabla$ is given (after a gauge transformation) as
\[\nabla=\begin{pmatrix} \nabla^S& \varphi^*+q\\ \varphi &  \nabla^{S^*}\end{pmatrix}.\]
 
\begin{remark}
We are mainly interested in the case that $q\neq0$. This case corresponds to non-trivial topological data $\mathcal C$, and the $SU(1,1)$-structure of $\nabla$ differs from that of 
$\nabla^{Fuchs}.$ In particular, $S^*\subset V$ is not the orthogonal complement of $S$ with respect to the indefinite hermitian metric induced by the flat  $\SU(1,1)$-connection $\nabla.$
\end{remark}
\subsection{Construction of real sections}\label{sec-con}
Consider the family of
gauge transformations (parametrised by $\lambda\in\C^*$)
\begin{equation}\label{eq:gga}g(\lambda):=\begin{pmatrix} 1&0\\0&\lambda\end{pmatrix}\end{equation}
with respect to the $C^\infty$-decomposition
$V=S\oplus S^*,$ and the family of flat connections
\[\nabla^\lambda:=\nabla.g(\lambda)=\begin{pmatrix} \nabla^S &0\\ 0&\nabla^{S^*} \end{pmatrix}+\lambda^{-1}\begin{pmatrix} 0&0\\ \varphi &0\end{pmatrix}+\lambda \begin{pmatrix} 0& \varphi^*+q\\0&0\end{pmatrix}.\]
We directly see that
 \begin{equation}\label{eq:graftingsection}\lambda\in\C\mapsto [\lambda,\bar\partial^{\nabla^\lambda},\lambda\partial^{\nabla^\lambda}]_\Sigma\end{equation}
defines a section $s$ of $\mathcal M_{DH}\to\C P^1$ over $\C\subset\C P^1.$
Over $\C^*\subset\C\subset\C P^1$  the section is also given by
\begin{equation}\label{constnab}\lambda\mapsto[\lambda,\bar\partial^\nabla,\lambda\partial^\nabla]_\Sigma.\end{equation}
\begin{lemma}\label{prorealsec}
Let $\nabla$ be a oper with $\SU(1,1)$-monodromy. Then, the section in \eqref{eq:graftingsection} extends holomorphically 
to $\lambda=\infty$. The extension is $\tau$-real.
\end{lemma}
\begin{proof}
We have that $\SU(1,1)$ and $\SL(2,\R)$ are conjugate.
Therefore, and because $\nabla$ has $\SU(1,1)$-monodromy, there exists a $\SL(2,\C)$-gauge transformation $h$ such that
\[\overline{\nabla}=\nabla.h.\]
Over  $\C^*$  the section is given by \eqref{constnab}.
Thus, we have
\begin{equation}
\begin{split}
\tau(s(\lambda))=&\tau([\lambda,\bar\partial^\nabla,\lambda\partial^\nabla]_\Sigma)\\
=&[-\bar\lambda,\overline{\bar\partial^\nabla},-\bar\lambda\overline{\partial^\nabla}]_{\bar\Sigma}\\
=&[-\bar\lambda,\partial^\nabla.h,-\bar\lambda{\bar\partial^\nabla}.h]_{\bar\Sigma}\\
=&[-\bar\lambda^{-1},\bar\partial^\nabla.h,-\bar\lambda^{-1}\partial^\nabla.h]_{\Sigma}\\
=&[-\bar\lambda^{-1},\bar\partial^\nabla,-\bar\lambda^{-1}\partial^\nabla]_{\Sigma}=s(-\bar\lambda^{-1})\\
\end{split}
\end{equation}
for all $\lambda\in\C^*.$
This shows that the  holomorphic section is $\tau$-real over $\C^*.$ As the section extends holomorphically to $\lambda=0,$ it follows from
$\tau$-reality that the section also extends holomorphically to $\lambda=\infty.$
\qed\end{proof}
Note that image of the section constructed in Lemma \ref{prorealsec} is contained in the smooth part of the Deligne-Hitchin moduli space.
\begin{definition}
A section of the form \eqref{eq:graftingsection} for an oper $\nabla$ with  real monodromy is called grafting section. It is called the uniformization section if it is given by the uniformization oper $\nabla^{Fuchs}.$
\end{definition}
\begin{theorem}
A grafting section is a $\tau$-real negative holomorphic section. It is admissible if and only if it is the uniformization section. Sections corresponding to different grafting data $\mathcal C$ and $\tilde{\mathcal C}$  are different.
\end{theorem}
\begin{proof}
We have already seen in Lemma \ref{prorealsec} that a grafting section $s$ is $\tau$-real. 
Due to the $\SU(1,1)$-structure there exists a $\lambda$-independent $\SL(2,\C)$-gauge transformation $h$ such that
\[\overline{\nabla}=\nabla.h.\]
We have $\bar h h=\pm\text{Id}$.
If $\bar h h=-\text{Id}$, $\nabla$ would be a flat $\SU(2)$ connection.
 Because $\nabla$ is an irreducible $\SU(1,1)$-connection we therefore have $\bar h h=\text{Id}$. Using
$\nabla^\lambda=\nabla.g(\lambda)$ for $g$ as in \eqref{eq:gga} we obtain
\[\overline{\nabla^{-\bar\lambda^{-1}}}=\nabla^\lambda .g^{-1}(\lambda)h\overline{g(-\bar\lambda^{-1})}.\]
Therefore \begin{equation}\label{eq:graftinggauge}\tilde g(\lambda):=ig^{-1}(\lambda)h\overline{g(-\bar\lambda^{-1})}\end{equation} satisfies
\[\overline{\tilde g(-\bar\lambda^{-1})}\tilde g(\lambda)=\text{Id}\]
(see also \cite[Lemma 2.15]{BHR} for the simple computation), and
\[\hat g(\lambda):=\lambda\tilde g(\lambda)\]
is a $\C^*$-family of $\SL(2,\C)$-valued gauge transformations such that
\[\overline{\nabla^{-\bar\lambda^{-1}}}=\nabla^\lambda .\hat g(\lambda)\]
and \[\overline{\hat g(-\bar\lambda^{-1})}\hat g(\lambda)=-\text{Id}.\]
Hence, $s$ is negative by definition. 

The uniformization section is given by a particular solution of the self-duality equations.
Hence, the uniformization section is admissible. Conversely, consider a grafting section $s$, and assume it is admissible. By \cite[Theorem 3.6]{BHR} the section corresponds to a solution of the self-duality equations. Therefore, it is a twistor line.
As the twistor line through a point in the Deligne-Hitchin moduli space is unique, and because the uniformization section and the grafting section coincide at $\lambda=0,$ $s$ must be the uniformization section.
Finally, if two sections are obtained by different graftings,  the monodromies of their flat  connections (at $\lambda=1$) are different. Therefore the sections are  different.
\qed\end{proof}

\begin{remark}
The extension of the $\C^*$-orbit  \ref{constnab} to $\lambda=0$ is a special instance of Simpson's construction \cite{Simpson-ite}. If the holomorphic bundle $\bar\partial^\nabla$ underlying a flat connection $\nabla$ with $\mathrm{SL}(2,\R)$-monodromy is stable, the $\C^*$-orbit yields a  $\tau$-real  positive section.
\end{remark}

\subsection{Real holomorphic sections and singular solutions of the self-duality equations}\label{subsec:rhs}
In \cite{BHR} it was shown that an admissible negative $\tau$-real holomorphic section is given by a solution of Hitchin's self-duality equations. The idea of proof is as follows: Take a lift of the section $s$ given by a $\C^*$ family of flat connections
\[\nabla^\lambda=\lambda^{-1}\Phi+\nabla+\dots\] where $(\bar\partial^\nabla,\Phi)$ is a stable Higgs pair. As $s$ is negative $\tau$-real, there exists a family of $\SL(2,\C)$ gauge transformations $g(\lambda)$ satisfying \eqref{eq:sign} and \eqref{eq:sign2} for the minus sign. Because $s$ is admissible, for every $p\in\Sigma$ the loop $\lambda\mapsto g_p(\lambda)$ is in the open big cell, hence admits a factorisation into a positive and a negative loop
(see \cite{Pressley_Segal_1986} for details about loop groups  or \cite{HH2} for a short summary of relevant material). Using a normalisation
(for example given by the reality condition), we have a splitting
\[g(\lambda)=g^+(\lambda)g^-(\lambda)\]
globally on $\Sigma$. 
It can be seen easily that
\[\nabla^\lambda.g^+(\lambda)\] is
the associated family of a solution to the self-duality equations (for an appropriate hermitian metric), see \cite[Theorem 3.6]{BHR} for details. 

If a $\tau$-real negative holomorphic section is not admissible, then (for a lift $\nabla^\lambda$ and the corresponding gauge satisfying \eqref{eq:sign} and \eqref{eq:sign2}) the loop $\lambda\mapsto g_p(\lambda)$ is not in the big cell for all $p\in\Sigma.$ Assume that $g_p$ is contained in the big cell
for all $p\in U$ in an open dense subset  $U\subset\Sigma$, and  intersects the first small cell transversely on $\Sigma\setminus U$. Then
the equivariant harmonic map to hyperbolic 3-space corresponding to the solution of the self-duality equations defined on $U$ (by the above procedure) extends smoothly through the boundary at infinity $S^2$ of $\mathbb H^3$, and yields a smooth map \[f\colon \Sigma\to \mathbb H^3\cup S^2\cup \mathbb H^3=S^3,\]
see \cite[$\S$ 5]{HH2}.
\begin{lemma}\label{lemsing}
Let $s$ be  a grafting section $s$. For $p\in\Sigma$ the corresponding loop $\tilde g_p$ in \eqref{eq:graftinggauge} is in the big cell  if $S_p$ and $S^\#_p$ do not coalesce, where $S^\#_p$ is the orthogonal bundle of $S$ with respect to the indefinite hermitian metric. This happens on an open dense subset $U\subset\Sigma$.
\end{lemma}
\begin{proof}
Recall the construction of the section of the Deligne-Hitchin moduli space in Section \ref{sec-con}.
Instead of taking the complementary line bundle $S^*=\ker\Phi$ of  $S$ provided by uniformization, we can as well use a different
complementary line bundle $L$ of $S.$ Define $g$ as in \eqref{eq:gga} but with respect to $S\oplus L.$
 This construction clearly yields the same section of the Deligne-Hitchin moduli space. The  two families of flat connections provided by two different choices of complementary line bundles
 differ by a family of $\SL(2,\C)$ gauge transformations which holomorphically extends to $\lambda=0.$ 
 
Consider a point $p\in\Sigma$ where $S_p$ and $S^\#_p$ do not coalesce. Consider a complementary line bundle $L$ of $S$ which coincides with $S^\#$ in a open neighbourhood $U$ of $p.$ A short computation
 then shows that  the corresponding map $\tilde g$ in \eqref{eq:graftinggauge} is  $\tilde g=\lambda^{-1}\delta$ for some $\delta\colon U\to \SL(2,\C).$ 
 Thus $\hat g_q(\lambda)=\lambda\tilde g_q(\lambda)$ is in the big cell for all $q\in U$.
 
Note that $S=S^\#$ cannot hold globally on $\Sigma$ because otherwise $S$ would be a parallel subbundle with respect to $\nabla$, which is a contradiction.
 The second statement of the lemma follows from the fact that $S$ and $S^\#$ are real analytic subbundles.
 \qed\end{proof}
This lemma shows that a grafting section gives rise to a solution of the self-duality equations on an open dense subset $U\subset\Sigma$. In particular, we obtain an equivariant harmonic map into $\mathbb H^3$ defined on the universal covering of $U$.

\begin{proposition}\label{devlemma}
The (equivariant) harmonic map with singularities given by a grafting section is given by the developing map of the corresponding real projective structure.
\end{proposition}
\begin{proof}
A grafting section has nilpotent Higgs field. Thus, the (equivariant) harmonic map corresponding to the solution of the self-duality equation  on the open dense subset $U\subset\Sigma$ provided by the previous lemma is a (equivariant) conformally parametrised minimal surface. In particular, it extends to a (equivariant) Willmore surface \cite[$\S5$]{HH2}.
Moreover, the Hopf differential of the surface does vanish, as a consequence of the fact that $\nabla^\lambda$ are gauge equivalent for all $\lambda\in\C^*$.
In \cite{BuPP} it is shown that for vanishing Hopf differential, the  developing map of the projective structure is equal to
 the (equivariant) Willmore immersion, and the proposition follows.
\qed\end{proof}
\begin{remark}\label{Rem:singset}
The previous proposition tells us that the loop $\tilde g$ in \eqref{eq:graftinggauge} is not contained in the big cell exactly at those points which are mapped to the boundary circle at infinity $S^1$ of $\mathbb H^2\subset\C P^1$ via the developing map of the real projective structure. Moreover, it can  be shown shown that $\tilde g$ in \eqref{eq:graftinggauge} intersects
the first small cell transversely, compare with \cite[Section 5]{HH2}.
\end{remark}
\subsection{Normal bundles of grafting sections}
The next goal is to compute the normal bundle of a grafting section.
We first recall the deformation theory of irreducible flat connections $\nabla$. The exterior differential $d^\nabla$
induces a complex
\[0\to\Omega^0(\Sigma,\text{End}_0(V))\to\Omega^1(\Sigma,\text{End}_0(V))\to\Omega^2(\Sigma,\text{End}_0(V))\to 0.\]
The tangent space at $\nabla$ of the deRham moduli space of flat irreducible connections  is naturally identified with
the first cohomology group
\[T_{[\nabla]}\mathcal M_{dR}=H^1(\Sigma,d^\nabla):=\text{ker} ({d^\nabla\colon\Omega^1\to\Omega^2})/\text{im}(d^\nabla\colon\Omega^0\to\Omega^1).\]
For an irreducible unitary flat connection the space $H^1(\Sigma,d^\nabla)$ can be identified with the space
of harmonic $\mathfrak{sl}(2,\C)$-valued 1-forms $\Psi,$ i.e., $d^\nabla\Psi=d^\nabla*\Psi=0$. The space of  harmonic 1-forms equals the direct sum of the space of Higgs fields with respect to $\bar\partial^\nabla$ and of
 the space of anti-Higgs fields with respect to $\partial^\nabla.$ In the case of a oper connection, the underlying holomorphic structure is unstable and admits  trace-free endomorphisms $0\neq X\in\Gamma(\Sigma,\text{End}_0(V))$
 with 
 \[\bar\partial X=0.\] By irreducibility of $\nabla$
 \[d^\nabla X=\partial^\nabla X\neq0\] is non-trivial, and by flatness of $\nabla$ it is holomorphic, i.e.,
 \begin{equation}\label{parX}\partial^\nabla X\in H^0(\Sigma, K\text{End}_0(V))\end{equation}
 is a Higgs field. Therefore, not every harmonic 1-form yields a non-trivial tangent direction of the deRham moduli space. The following lemma shows that for an oper connection $\nabla$ with real monodromy every
 tangent vector \[\Psi\in T_{[\nabla]}\mathcal M_{dR}=H^1(\Sigma,d^\nabla)\]
 can be represented by special harmonic 1-forms, i.e., the direct sum of a nilpotent Higgs field with a nilpotent anti-Higgs field.
Denote by
 \[\mathcal Q:=\{\Phi\in H^0(\Sigma, K\text{End}_0(V),\bar\partial^\nabla)\mid \det \Phi=0\}\]
 and
 \[\bar{\mathcal Q}:=\{\Psi\in H^0(\bar\Sigma, K_{\bar\Sigma}\text{End}_0(V),\partial^\nabla)\mid \det \Psi=0\}\]
 the space of nilpotent Higgs fields respectively anti-Higgs fields.
 Recall the form \eqref{nablafudec}. Thus $\bar\partial^\nabla$ is upper triangular with respect to $V=S\oplus S^*$.
 A direct computation using that $\varphi^*$ gives a nontrivial  element in $H^1(\Sigma,K^{-1})$ shows that every Higgs field of $(V,\bar\partial^\nabla)$ must be upper 
 triangular as well. Thus, $\mathcal Q$ consists of upper triangular Higgs fields with respect to $V=S\oplus S^*.$
 In particular  $\mathcal Q$ is a vector space modelled on the space of holomorphic quadratic differentials, and its dimension is $3g-3,$  where $g$ is the genus of $\Sigma$. Analogusly, $\bar{\mathcal Q}$  is also a vector space of dimension $3g-3$.

\begin{lemma}\label{lem:ts}
Let $\nabla$ be an oper  with real monodromy. Then the natural map
\[\mathcal Q\oplus\bar {\mathcal Q}\to H^1(\Sigma,d^\nabla)=\text{ker} ({d^\nabla\colon\Omega^1\to\Omega^2})/\text{im}(d^\nabla\colon\Omega^0\to\Omega^1)\]
is an isomorphism.
\end{lemma}
\begin{proof}
Higgs fields and anti-Higgs fields are clearly closed with respect to $d^\nabla,$ so we have a well-defined map
$\mathcal Q\oplus\bar {\mathcal Q}\to H^1(\Sigma,d^\nabla)$.
As the dimension is $6g-6$ for both spaces it remains to show that this map is injective.
As elements of $\mathcal Q$ are $(1,0)$-forms on $X$ and  elements of $\bar{\mathcal Q}$ are $(0,1)$-forms on $X$
the  map $\mathcal Q\oplus\bar {\mathcal Q}\to\text{ker} ({d^\nabla\colon\Omega^1\to\Omega^2})$
is injective.
Let \[\Psi\in \mathcal Q\oplus\bar {\mathcal Q}\] be in the kernel of the map
$\mathcal Q\oplus\bar {\mathcal Q}\to H^1(\Sigma,d^\nabla)$, i.e., there exists
\[X\in\Gamma(\Sigma, \text{End}_0(V))\]
with
\[\Psi=d^\nabla X.\]
We claim that $\Psi=0$. For any endomorphism $A$ and any endomorphism-valued 1-form $\omega$ we denote
by $A^\#$ and $\omega^\#$ the adjoint endomorphism and the adjoint 1-form with respect to the indefinite hermitian metric.
We  decompose
\[\Psi=\Psi^++\Psi^-\]
into its  hermitian symmetric part $\Psi^+=\tfrac{1}{2}(\Psi+\Psi^\#)$ and skew hermitian part $\Psi^-=\tfrac{1}{2}(\Psi-\Psi^\#)$ 
with respect to the indefinite hermitian metric, and analogously 
\[X=\tfrac{1}{2}(X+X^{\#})+\tfrac{1}{2}(X-X^{\#})=:X^++X^-.\]

Note that for $\Psi\in \mathcal Q\oplus\bar {\mathcal Q}$ also $\Psi^\#\in \mathcal Q\oplus\bar {\mathcal Q}.$
 Moreover, since $\nabla$ is unitary with respect to the 
 indefinite hermitian metric,
 we obtain
\[d^\nabla (X^\#)=(d^\nabla X)^\#=\Psi^\#\in\mathcal Q\oplus\bar {\mathcal Q}.\]
This gives  \[d^\nabla X^\pm =\Psi^\pm \in\mathcal Q\oplus\bar {\mathcal Q},\] and 
it suffices  to prove the claim for  $\Psi^\pm\in\mathcal Q\oplus\bar{\mathcal Q}$.  

We 
consider $\Psi^-$ only, the proof for $\Psi^+$ works completely analogously.
Consider the Higgs field
\[\Phi:=(\Psi^-)^{(1,0)}=\partial^\nabla X^-\in\mathcal Q.\]
Then we have
\[\Phi^\#=-(\Psi^-)^{(0,1)}=-\bar\partial^\nabla X^-\in\bar{\mathcal Q},\]
and integration yields
\[\int_\Sigma\text{tr}(\Phi\wedge\Phi^\#)=-\int_\Sigma\text{tr}(\partial^\nabla X^-\wedge\bar\partial^\nabla X^-)=-
\int_\Sigma d\text{tr}(X^-\bar\partial^\nabla X^-)=0.\]
For a non-degenerate indefinite hermitian metric on a vector space of dimension 2 every non-vanishing nilpotent endomorphism $A$ satisfies
\begin{equation}\label{eq:nilsign}\text{tr}(AA^\#)\geq 0\end{equation}
with equality if and only if the kernel of $A$ is  a null-line. By construction of the oper $\nabla$ with real monodromy and by Remark \ref{Rem:singset}, the fibers of the kernel bundle of $\Phi$ (given by the holomorphic line bundle $S$) are null exactly 
where the developing map of the real projective structure crosses
the boundary of the hyperbolic disc. This is a subset of measure $0$, therefore
\[\int_\Sigma\text{tr}(\Phi\wedge\Phi^\#)=0\] holds
if and only if $\Phi=0.$
Thus,  $\Psi^-=0$ and analogously  $\Psi^+=0$ proving the lemma.
\qed\end{proof}

In the $\SL(2,\C)$-case, 
the normal bundle of a holomorphic section $s\colon \C P^1\to\mathcal M_{DH}$ is a holomorphic vector bundle
over $\C P^1$
of rank $6g-6$.  
For the uniformization section (as for all twistor lines), the normal bundle is 
\[\mathcal O(1)^{6g-6}\to\C P^1,\]
which is fundamental for the twistor approach to hyperK\"ahler manifolds \cite{HKLR}. 
\begin{theorem}\label{the:normal}
The normal bundle $\mathcal N$ of a grafting section on a Riemann surface $\Sigma$ of genus $g$ is
\[\mathcal O(1)^{6g-6}\to\C P^1.\]
\end{theorem}
\begin{proof}
We first construct a (holomorphic) bundle homomorphism
\[\mathcal H\colon\mathcal O(1)\otimes (\mathcal Q\oplus \bar{\mathcal Q})\to\mathcal N.\]
Then we show that $\mathcal H$ is an isomorphism of holomorphic bundles.

Consider \[0\neq Q\in H^0(\Sigma, K^2),\]
which determines a non-zero element in $\mathcal Q$.
Over $\C$ the grafting section is given by
\[s(\lambda)=[\lambda,\begin{pmatrix} \bar\partial^S&\lambda\varphi^*\\0&\bar\partial^{S^*}\end{pmatrix},\lambda\begin{pmatrix} \partial^S&\lambda q \\ \lambda^{-1}\varphi&\partial^{S^*}\end{pmatrix}]_\Sigma.\]
Consider the infinitesimal deformation
 \[s^Q_t(\lambda)=[\lambda,\begin{pmatrix} \bar\partial^S&\lambda\varphi^*\\0&\bar\partial^{S^*}\end{pmatrix},\lambda\begin{pmatrix} \partial^S&\lambda q \\ \lambda^{-1}\varphi&\partial^{S^*}\end{pmatrix}+t\begin{pmatrix} 0 &Q\\0&0\end{pmatrix}]_\Sigma.\]
 By Lemma \eqref{lem:ts} this gives a non-vanishing holomorphic section $s_Q$ of the normal bundle of $s$ over $\C$.
 Note that  over $\C^*$ this normal field is given by
 \begin{equation}\label{normalq2}s^Q_t(\lambda)=[\lambda,\begin{pmatrix} \bar\partial^S&\varphi^*\\0&\bar\partial^{S^*}\end{pmatrix},\lambda\begin{pmatrix} \partial^S& q \\ \varphi&\partial^{S^*}\end{pmatrix}+t\lambda^{-1}\begin{pmatrix} 0 &Q\\0&0\end{pmatrix}]_\Sigma.\end{equation}
We want to analyse the behaviour of this normal bundle section at $\lambda=\infty.$ 

Likewise, we can start with an element 
\[\bar Q=\begin{pmatrix} x&y\\z&-x \end{pmatrix}\in\bar{\mathcal Q},\]
where $x\in\Gamma(\Sigma,\bar K),$ $y\in\Gamma(\Sigma,\bar K K)$ and $z\in\Gamma(\Sigma,\bar KK^{-1})$. We do the analogous construction to obtain a holomorphic normal field of $s$ over $\C P^1\setminus\{0\},$ and want to study its behaviour at $\lambda=0.$ 
By comparing with \eqref{normalq2}, 
this normal field is over $\C^*$  given by
\[s^{\bar Q}_t(\lambda)=[\lambda,\begin{pmatrix} \bar\partial^S&\varphi^*\\0&\bar\partial^{S^*}\end{pmatrix}+t\lambda^2\begin{pmatrix} x&y\\z&-x\end{pmatrix},\lambda\begin{pmatrix} \partial^S& q \\ \varphi&\partial^{S^*}\end{pmatrix}]_\Sigma.\] Hence, over
$\C$  the normal section is given by
\begin{equation}\label{eq:norsec}s^{\bar Q}_t(\lambda)=[\lambda,\begin{pmatrix} \bar\partial^S&\lambda \varphi^*\\0&\bar\partial^{S^*}\end{pmatrix}+t\lambda^2\begin{pmatrix} x&\lambda y\\ \lambda ^{-1}z&-x\end{pmatrix},\lambda\begin{pmatrix} \partial^S& \lambda q \\ \lambda^{-1}\varphi&\partial^{S^*}\end{pmatrix}]_\Sigma.\end{equation}
From \eqref{eq:norsec}, the normal bundle section $s_{\bar Q}$ given by the variation $s^{\bar Q}_t$ extends holomorphically to $\lambda=0$ with a 
zero of order at least 1. Analogously, or by using the $\SU(1,1)$-structure, we see that the normal field obtained by  $s^{Q}_t$ gives a section 
which extends holomorphically to $\infty$ with a zero of order at least 1.  We obtain a $6g-6$ 
dimensional space (parametrised by $\mathcal Q\oplus\bar{\mathcal Q}$) of holomorphic sections.
Because of Lemma \ref{lem:ts}, the evaluation at $\lambda$ spans $\mathcal N_\lambda$ for all $\lambda\in\C^*.$
We claim that the section given by elements in $\mathcal Q$ respectively $\bar{\mathcal Q}$  have first order zeros at $\lambda=\infty$ respectively $\lambda=0$.
 In fact, this follows  by another application of Lemma \ref{lem:ts}: The space of  Higgs fields $\Phi$ is the direct sum of $\mathcal Q$ and the image of the space of
holomorphic trace-free endomorphisms under the map $\partial^\nabla$, as can be deduced from Riemann-Roch and \eqref{parX}.
Assume \[ \Psi=\begin{pmatrix}x&y\\z&-x\end{pmatrix}\in\bar{\mathcal Q}\] satisfies
\[\int_\Sigma Qz=0\] for all $Q\in H^0(\Sigma,K^2)$. Then, $\int\text{tr}(\Psi\wedge \Phi)=0$ for every holomorphic Higgs field $\Phi,$
and
Serre-duality implies that 
$\Psi$ is in the image of $\bar\partial^\nabla.$ By  Lemma \ref{lem:ts} $\Psi=0.$

We have a well-defined holomorphic bundle homomorphism via
\[\mathcal H\colon\mathcal O(1)\otimes (\mathcal Q\oplus \bar{\mathcal Q})\to\mathcal N,\quad (a+b\lambda)(Q\oplus\bar Q)\mapsto
(a+b\lambda)s_Q+ (a\lambda^{-1}+b) s_{\bar Q}.\]
The map $\mathcal H$ is an isomorphism over $\C^*$.
We have already deduced from Serre duality  that the pairing
\[(Q,\bar Q)\mapsto \int_\Sigma Qz\]
is a duality, where $z$ is the lower left entry (with respect to $S\oplus S^*$) of the nilpotent anti-Higgs field determined by $\bar Q.$
Hence, $\mathcal H$ is an isomorphism at $\lambda=0$ as well. Finally, the $\SU(1,1)$-symmetry shows that $\mathcal H$ is an isomorphism at $\lambda=\infty$ as well.
 \qed\end{proof}
\section{Hyper-K\"ahler components}\label{sec:hk}
We are now able to state and prove our main theorem:
\begin{theorem}\label{MT}
For a real projective structure $P$ on a Riemann surface there exists a $6g-6$-dimensional manifold $\mathcal M_P$ of real holomorphic sections
of the Deligne-Hitchin moduli space. Moreover, $\mathcal M_P$  is equipped with a hyper-K\"ahler metric.

There exists a neighbourhood $\mathcal U$ of $s_P\in\mathcal M_P$, such that every section $s\in\mathcal U$ gives rise to a solution of the 
self-duality equation away from a singularity set $S$. The set $S$ is given by a collection of smooth simple curves in $\Sigma$ which are homotopic to the singular curves 
of the constant curvature -1 metric associated to $P.$
\end{theorem}
\begin{proof}
The first part of the theorem  is an application of \cite[Theorem 3.3]{HKLR} using Theorem \ref{the:normal}:  
it remains to show that the induced hyper-K\"ahler metric is positive definite. 
Since our construction is local it is sufficient to compute the induced bilinear form at a grafting section.
By Theorem \ref{the:normal} and its proof, a complex tangent vector $X$ to the space of holomorphic sections at the grafting section
is given by
\begin{equation}\label{tangentenvektor}X=(P_0+P_1\lambda,Q_0+Q_1\lambda)\end{equation}
with
\[P_0,P_1\in \bar{\mathcal Q} \quad \text{and}\quad Q_0,Q_1\in {\mathcal Q}.\]
As in \cite{HKLR}, the twisted symplectic form induces a complex bilinear  form $g_\C$
on the space of holomorphic sections. In the case of the Deligne-Hitchin moduli space and $X$ as in \eqref{tangentenvektor}, it is given by
\[g_\C(X,X)=2i\int_\Sigma\text{tr}(-Q_1\wedge P_0+Q_0\wedge P_1),\]
see  \cite[Equation 6.2]{HitchinSD} and also \cite[$\S 2$ and $\S4$]{GeometrySections}.
Using the definition of $\tau$ in \eqref{deftaus}, a tangent vector $X$ as in \eqref{tangentenvektor} is real, i.e. a tangent vector to the space of real holomorphic sections,
if and only 
\[Q_1=-P_0^{\#} \quad\text{and}\quad P_1=Q_0^{\#}.\]
Using \eqref{eq:nilsign} 
we get  analogous to the proof of Lemma \ref{lem:ts} 
\[g_\C(X,X)\geq 0 \]
for real tangent vectors $X$ with equality if and only if $X=0.$
Therefore, the induced bilinear form on the space of real sections is positive definit at a grafting section.
By continuity, this also holds on a suitable local neighbourhood -- denoted by $\mathcal M_P$ -- of the grafting section in the space of real holomorphic sections.
 By \cite[Theorem 3.3]{HKLR}, we obtain a (Riemannian) hyper-K\"ahler metric on $\mathcal M_P$.
\\

The second part of the theorem follows from the discussion in Section \ref{subsec:rhs} together with the following facts:
\begin{itemize}
\item the choice of lifts of sections  $s$ can be done in a smooth way (depending on $s$) locally;
\item by irreducibility, the gauges satisfying  \eqref{eq:sign} depend smoothly on $s$;
\item being in the big cell (respectively being in the union of the big cell and the first small cell) is an open condition.
\end{itemize}
The singularity set of the solution of the self-duality equation is given by the set of points $p\in\Sigma$ where the
gauge satisfying  \eqref{eq:sign} is not in the big cell. By the above bullet points and because of Lemma \ref{lemsing} and Proposition \ref{devlemma}, this happens along curves which are homotopic
to the singular curves of the constant curvature -1 metric provided by the real projective structure $P$.
\qed\end{proof}
\begin{remark}\label{lastremark}
We cannot prove that the components of  real holomorphic sections associated to different real projective structures 
are different. If the singularity locus of the self-duality solutions associated to a real holomorphic section would be an invariant of its component of real holomorphic sections
these spaces would actually be different.
\end{remark}
In \cite{BeHR}, an energy functional $\mathcal E$ on the space of holomorphic sections of the Deligne-Hitchin moduli space has been defined.
It takes real values on real holomorphic sections. Moreover, by \cite[Theorem 3.13]{BeHR} this energy is a K\"ahler potential for the hyper-K\"ahler metric constructed in Theorem \ref{MT}. For a grafting section $s$, $\mathcal E(s)=1-g$ where $g$ is the genus of the 
surface.

\bibliographystyle{amsplain}
\bibliography{references}
\end{document}